\documentclass[10pt]{article}
\usepackage{amsfonts}
\usepackage{bbm}
\usepackage{latexsym,bm}
\usepackage{amssymb}
\usepackage{amsmath}
\usepackage{amsthm}
\usepackage[all]{xy}
\theoremstyle{plain}
\newtheorem{thm}{Theorem}[section]

\newtheorem{lem}[thm]{Lemma}
\newtheorem{prop}[thm]{Proposition}
\newtheorem{cor}[thm]{Corollary}
\theoremstyle{remark}
\newtheorem{e.g.}[thm]{Example}
\newtheorem{rmk}[thm]{Remark}
\theoremstyle{definition}

\begin{document}

\title{\Large Subvarieties of small codimension in smooth projective varieties}
\author{\large Qifeng Li}
\maketitle

\begin{abstract}
Let $X\subsetneq\mathbb{P}_{\mathbb{C}}^{N}$ be an $n$-dimensional nondegenerate smooth projective variety containing an $m$-dimensional subvariety $Y$. Assume that either $m>\frac{n}{2}$ and $X$ is a complete intersection or that $m\geq\frac{N}{2}$, we show $\text{deg}(X)|\text{deg}(Y)$ and $\text{codim}_{\langle Y\rangle}Y\geq\text{codim}_{\mathbb{P}^{N}}X$, where $\langle Y\rangle$ is the linear span of $Y$. These bounds are sharp. As an application, we classify smooth projective $n$-dimensional quadratic varieties swept out by $m\geq[\frac{n}{2}]+1$ dimensional quadrics passing through one point.
\end{abstract}


\setcounter{section}{0}

\section{\large Introduction}

We work over the complex number field $\mathbb{C}$. Varieties in this paper are assumed to be irreducible, unless
otherwise stated. Let $X\subsetneq\mathbb{P}^{N}$ be an $n$-dimensional  nondegenerate smooth projective variety. For a closed subvariety $Y\subseteq X$, we denote by $\langle Y\rangle$ its linear span. It turns out that in general $Y$ has not bigger codimension in $\langle Y\rangle$ than $X$ does in $\mathbb{P}^{N}$. We shall show that $\text{codim}_{\langle Y\rangle}Y\geq\text{codim}_{\mathbb{P}^{N}}X$ if $\text{dim}(Y)$ is big enough.
More precisely, we have the following:

\begin{thm} \label{thm: introduction deg(X)|deg(Y)}
Let $X\subsetneq\mathbb{P}^{N}$ be a nondegenerate smooth projective variety of dimension $n$ and let $Y\subseteq X$ be a subvariety of dimension $m$. Assume that either

$(i)$ $X$ is a complete intersection in $\mathbb{P}^{N}$, and $m>\frac{n}{2}$; or

$(ii)$ $m\geq \frac{N}{2}$.

Then $\text{deg}(X)|\text{deg}(Y)$, and $\text{codim}_{\langle Y\rangle}Y\geq\text{codim}_{\mathbb{P}^{N}}X$.
\end{thm}

These bounds are sharp as shown by complete intersections of two hyperquadrics  (see Section \ref{section: complete intersection of two hyperquadrics}).

A subvariety  $X\subseteq\mathbb{P}^{N}$ is called \emph{quadratic} if $X$ is the scheme theoretic intersection of quadric hypersurfaces.
We say that $X$ is swept out by $m$-dimensional quadrics passing through a point if there is a fixed point $x_{0}\in X$ such that for a general point
$x\in X$, there exists an $m$-dimensional quadric $Q_{x}$ satisfying
$x_{0}, x\in Q_{x}\subseteq X$.  Recall that Hartshorne Conjecture (\cite{Ha}) says that if $X\subseteq\mathbb{P}^{N}$ is a smooth projective variety of dimension $n>\frac{2N}{3}$, then $X\subseteq\mathbb{P}^{N}$ is a complete intersection. If $X$ is quadratic, then this conjecture is proved by Ionescu and Russo in \cite{IR}.

As an application of Theorem \ref{thm: introduction deg(X)|deg(Y)}, we obtain the classification of quadratic varieties swept out by quadric hypersurfaces of dimension $[\frac{n}{2}]+1$, providing a slight improvement of \cite[Prop. 3]{Fu}.

\begin{thm} \label{thm: introduction application to var. swept out by large dim. quadrics}
Let $X\subsetneq\mathbb{P}^{N}$ be an $n$-dimensional nondegenerate smooth
projective subvariety swept out by quadric hypersurfaces of dimension $[\frac{n}{2}]+1$ passing though a fixed point of $X$.  Assume that either Hartshorne Conjecture on Complete Intersection holds or that $X \subseteq \mathbb{P}^N$ is a quadratic variety, then $X \subseteq \mathbb{P}^N$ is projectively equivalent to one of the following :

$(a)$ a quadric hypersurface;

$(b)$ the Segre 3-fold
$\mathbb{P}^{1}\times\mathbb{P}^{2}\subseteq\mathbb{P}^{5}$;

$(c)$ the Pl\"{u}cker embedding
$\mathbb{G}(1,4)\subseteq\mathbb{P}^{9}$;

$(d)$ the 10-dimensional spinor variety
$S^{10}\subseteq\mathbb{P}^{15}$;

$(e)$ a general hyperplane section of $(b)$ or $(c)$.
\end{thm}

By Proposition \ref{prop: complete intersection of two hyperquadrics}, the bound $[\frac{n}{2}]+1$ in Theorem \ref{thm: introduction application to var. swept out by large dim. quadrics} is sharp.

The paper is organized as follows. In Section \ref{section: main theorms and their corollaries}, we prove Theorem \ref{thm: introduction deg(X)|deg(Y)} and its corollaries. In Section \ref{section: complete intersection of two hyperquadrics}, we study complete intersections of two hyperquadrics, which shows that the bounds in Theorem \ref{thm: introduction deg(X)|deg(Y)} are sharp. In Section \ref{section: applications to quadratic varieties}, we prove Theorem \ref{thm: introduction application to var. swept out by large dim. quadrics}.

\bigskip

\textbf{\normalsize Acknowledgements.} It's my great pleasure to thank professor Baohua Fu for introducing
this problem and giving me many useful suggestions. I want to thank professor Fedor L. Zak for useful discussions and it should be pointed out that the idea of the proof of Theorem \ref{thm: introduction deg(X)|deg(Y)} is due to him. I also want to
thank professor Francesco Russo for his encouragement. I am grateful to the reviewers e for helpful comments and for the careful reading.

\section{\large Subvarieties of small codimension} \label{section: main theorms and their corollaries}

For a variety $X$ in a projective space $\mathbb{P}^{N}$, denote by $\langle X\rangle$ the linear span of $X$ in $\mathbb{P}^{N}$. Let $Y, Z\subseteq\mathbb{P}^{N}$ be two closed subschemes. Let $I(Y)$ and $I(Z)$ be the homogeneous ideal of $Y$ and $Z$ respectively.
Denote by $Y\cap_{sch} Z$  (resp. $Y\cap Z$) the closed subscheme of $\mathbb{P}^{N}$ defined by the homogeneous ideal $I(Y)+I(Z)$ (resp. the radical homogeneous ideal $\sqrt{I(Y)+I(Z)}$, and we call it the scheme-theoretic intersection (resp. the set-theoretic intersection) of $Y$ and $Z$.

\begin{proof}[Proof of Theorem \ref{thm: introduction deg(X)|deg(Y)}]
$(i)$ By the Lefschetz theorem (see, for example, \cite[Thm. 2.1(a)]{Ha}), the natural morphism $\phi: H^{2n-2m}(\mathbb{P}^{N}, \mathbb{Z})\rightarrow H^{2n-2m}(X, \mathbb{Z})$ is an isomorphism. Let $\tilde{L}$ be a general $(N+m-n)$-dimensional linear subspace. Then $[Y]=r\phi([\tilde{L}])\in H^{2n-2m}(X, \mathbb{Z})$ for some positive integer $r$. Since the degree is a homological invariant, $\text{deg}(Y)=r \text{deg}(X)$.

\medskip

Assume $\text{dim}(L)< N-n+m$, where $L=\langle Y\rangle$. Take $Y'$ to be a subvariety of $X$ of maximal dimension satisfying $\text{dim}(L')< N-n+\text{dim}(Y')$, where $L'=\langle Y'\rangle$. Let $Y''=L'\cap X$ and $m''=\text{dim}(Y'')$. Then $\text{dim}(L')< N-n+m''$ and $\frac{n}{2}<m\leq m''\leq n-1$. Let $c=n-m''$. Take any $c$ hyperplanes $H_{1},\ldots, H_{c}$ containing $L'$ such that $\text{dim}(\bigcap\limits_{i=1}^{c}H_{i})=N-c$. Denote by $X_{0}=X$, $L_{0}=\mathbb{P}^{N}$, $L_{s+1}=L_{s}\cap H_{s+1}$, and $X_{s+1}=X_{s}\cap H_{s+1}$ for each $0\leq s\leq c-1$. By induction on $s$ and the choice of $Y'$, for each $0\leq s\leq c-1$,  $X_{s}$ is of pure dimension $n-s$ and each component of it is nondegenerate in $L_{s}$. Then $X_{c}=X_{c-1}\cap H_{c}=X\cap L_{c}$ is of pure dimension $m''$. We can choose $H_{1},\ldots,H_{c}$ satisfying that $L_{c}\cap X\supsetneq L'\cap X$, i.e. $X_{c}\supsetneq Y''$. Hence, $X_{c}$ is reducible. Denote by $X_{c}=\bigcup\limits_{i=1}^{t}W_{i}$, where these $W_{i}$ are the irreducible components of $X_{c}$. By the Lefschetz theorem, $\psi: H^{2c}(\mathbb{P}^{N}, \mathbb{Z})\rightarrow H^{2c}(X, \mathbb{Z})$ is an isomorphism. Thus, there is a positive integer $r_{i}$ such that $[W_{i}]=r_{i}\psi([L_{c}])\in H^{2c}(X, \mathbb{Z})$ for each $i$. In particular, $\text{deg}(W_{i})=r_{i}\text{deg}(L_{c}\cap_{sch}X)$. On the other hand, $\text{deg}(L_{c}\cap_{sch}X)\geq \text{deg}(X_{c})=\sum\limits_{i=1}^{t}r_{i}\text{deg}(W_{i})$ and $t\geq 2$, which is a contradiction.

$(ii)$ By the Barth-Larsen theorem in \cite{BL}, the natural morphisms $H^{i}(\mathbb{P}^{N}, \mathbb{Z})\rightarrow H^{i}(X, \mathbb{Z})$ are isomorphisms for $i\leq 2n-N$. An argument similar to  $(i)$ proves the conclusion.
\end{proof}

\begin{rmk}
It should be noticed that there was similar results due to Zak by another approach. Let $X\subsetneq\mathbb{P}^{N}$ be a nondegenerate smooth projective variety of dimension $n$ and let $Y\subseteq X$ be a subvariety of dimension $m$. Denote by $L=\langle Y\rangle$ and $r=\dim(L)$. Then by \cite[Prop. I. 2.26]{Za}, $r<\min{n-1, \frac{n-1}{2}}$.
\end{rmk}

\begin{cor} \label{cor: L cap X irrecucible and L not contained in X}
Let $X\subsetneq\mathbb{P}^{N}$ be a nondegenerate smooth projective variety of dimension $n$ and $L\subseteq\mathbb{P}^{N}$ be a linear subspace of dimension $r$. Assume that either

$(i)$ $X$ is a complete intersection in $\mathbb{P}^{N}$, and $r>N-\frac{n}{2}$; or

$(ii)$ $r\geq \frac{3N}{2}-n$.

Then $L\cap_{sch} X$ is irreducible, it has multiplicity one on the unique irreducible component, and
$\text{dim}(L\cap X)=r+n-N$. In particular, $L\nsubseteq X$.
\end{cor}

\begin{proof}
$(i)$ Denote by $L\cap X=\bigcup\limits_{i}^{t}W_{i}$, where $W_{i}$ are the irreducible components of $L\cap X$. We have $\text{dim}(W_{i})\geq r+n-N > \frac{n}{2}$. By Theorem \ref{thm: introduction deg(X)|deg(Y)}$(i)$, $N-n+\text{dim}(W_{i})\leq\text{dim}\langle W_{i}\rangle\leq r$. Thus, $\text{dim}(W_{i})=r+n-N$, i.e. $L\cap X$ is of pure dimension $r+n-N$. Then $\text{deg}(X)=\text{deg}(L\cap_{sch}X)\geq\sum\limits_{i=1}^{t}\text{deg}(W_{i})\geq t\, \text{deg}(X)$, where the last inequality follows from Theorem \ref{thm: introduction deg(X)|deg(Y)}. Hence, $t=1$ and $\text{deg}(W_{1})=\text{deg}(L\cap_{sch}X)$, i.e. $L\cap X$ is irreducible and has multiplicity one on the unique irreducible component.

$(ii)$ By Theorem \ref{thm: introduction deg(X)|deg(Y)}$(ii)$, an argument similar to the proof of $(i)$ proves the conclusion.
\end{proof}

\begin{rmk}
The bounds $r>N-\frac{n}{2}$ and $r\geq \frac{3N}{2}-n$ in Corollary \ref{cor: L cap X irrecucible and L not contained in X} are not optimal to guarantee $L\nsubseteq X$. In fact, if $X$ is an $n$-dimensional projective manifold which is not a projective space, then $X$ does not contain linear subspaces of dimension greater than $[\frac{N-1}{2}]$ (see \cite[Cor. I. 2.10]{Za}).
\end{rmk}

Now we want to study the cases with $\text{deg}(X)=\text{deg}(Y)$, where $Y$ is a subvariety of $X$ of small codimension. A projective variety $X\subseteq\mathbb{P}^{N}$ is called a projective hypersurface if $X$ is a hypersurface in $\langle X
\rangle$.

\begin{cor} \label{cor: hypersurface deg(X)=deg(Y)}
Let $X\subseteq\mathbb{P}^{n+1}$ be a nondegenerate smooth hypersurface. Assume that $X$ contains an $m>\frac{n}{2}$
dimensional projective hypersurface $Y$, then $\text{deg}(X)=\text{deg}(Y)$.
\end{cor}

\begin{proof}
We can assume that $Y\neq X$. Let $L=\langle Y\rangle$. By Corollary \ref{cor: L cap X irrecucible and L not contained in X}, $L\nsubseteq X$ and $L\cap X=Y$. Hence, $\text{deg}(Y)\leq \text{deg}(L\cap_{sch}X)=\text{deg}(X)$. By Theorem \ref{thm: introduction deg(X)|deg(Y)}, $\text{deg}(X)\leq\text{deg}(Y)$. Hence, $\text{deg}(X)=\text{deg}(Y)$.
\end{proof}

\section{\large Smooth complete intersection of two hyperquadrics} \label{section: complete intersection of two hyperquadrics}

The main result of this section is the following:

\begin{prop} \label{prop: complete intersection of two hyperquadrics}
Let $X\subseteq\mathbb{P}^{n+2}$ be a smooth complete intersection
of two hyperquadrics $X_{1}$ and $X_{2}$.

$(i)$ $X$ does not contain any $m=[\frac{n}{2}]+1$ dimensional
quadric.

$(ii)$ If $n$ is even and $m=\frac{n}{2}$, then for a general point
$x\in X$, there is an $m$-dimensional quadric $Q_{x}$ satisfying
that $x\in Q_{x}\subseteq X$.

$(iii)$ Assume $n\geq 3$. Let $x_{0}\in X$ be an arbitrary fixed point and $m=[\frac{n-1}{2}]$. Then for a general point $x\in X$, there is an
$m$-dimensional quadric $Q_{x}$ satisfying that $x_{0}, x\in
Q_{x}\subseteq X$.
\end{prop}

\begin{rmk} \label{rmk: bounds are sharp}
From Proposition \ref{prop: complete intersection of two hyperquadrics}, we see that the bound in Theorem \ref{thm: introduction deg(X)|deg(Y)}$(i)$ is sharp, and when $N$ is even, the bound in Theorem \ref{thm: introduction deg(X)|deg(Y)}$(ii)$ is also sharp. On the other hand, consider the Grassmann variety $\mathbb{G}(1, 4)\subsetneq\mathbb{P}^{9}$ parametering lines in $\mathbb{P}^{4}$. Let $H$ be a hyperplane in $\mathbb{P}^{4}$. Denote by $X=\mathbb{G}(1, 4)$ and $Y=\mathbb{G}(1, H)$. Then $\text{dim}(X)=6$, $\text{dim}(Y)=4$, $\text{deg}(X)=5$ and $\text{deg}(Y)=2$. Moreover, $X$ is nondegenerate in $\mathbb{P}^{9}$. This shows that when $N$ is odd, the bound in Theorem \ref{thm: introduction deg(X)|deg(Y)}$(ii)$ is also sharp.
\end{rmk}

For a variety $X$, denote by $X_{sm}$ and $\text{Sing}(X)$ the smooth locus and the singular locus respectively.

\begin{lem} \label{lem: X cap Sing(Z) contained in Sing(X)}
Assume that $X\subseteq\mathbb{P}^{N}$ is an $n$-dimensional
complete intersection of hypersurfaces $X_{1},\ldots,X_{c}$, where
$c=\text{codim}_{\mathbb{P}^{N}}(X)$. Let
$i_{1},\ldots,i_{\nu}\in\{1,\ldots,c\}$ be any $\nu$ distinct
numbers and let $Z$ be the complete intersection of
$X_{i_{1}},\ldots,X_{i_{\nu}}$. Then $X_{sm}\subseteq Z_{sm}$, i.e.
$X\cap \text{Sing}(Z)\subseteq \text{Sing}(X)$.
\end{lem}

\begin{proof}
Let $S=k[t_{0},\ldots,t_{N}]$ be the homogeneous
coordinate ring of $\mathbb{P}^{N}$. Then $I(X_{i})=(f_{i})$ and
$I(X)=(f_{1},\ldots,f_{c})$. At any regular closed point $x\in
X_{sm}$, the rank of the matrix $(\frac{\partial f_{i}}{\partial
t_{j}})_{c\times (N+1)}$ is $c$, hence the rank of the matrix
$(\frac{\partial f_{i_{l}}}{\partial t_{j}})_{\nu\times (N+1)}$ is
$\nu$. Therefore, $x\in Z_{sm}$.
\end{proof}

\begin{lem} \label{lem: general Z contain smooth X is smooth}
Assume that $X\subseteq\mathbb{P}^{N}$ is an $n$-dimensional smooth complete
intersection of hypersurfaces $X_{1},\ldots,X_{c}$, where
$c=\text{codim}_{\mathbb{P}^{N}}(X)$ and $\text{deg}(X_{1})=\cdots=\text{deg}(X_{c})=d$.
Let $Z_{1},\ldots,Z_{t}$ be general $t\leq c$ hypersurfaces of
degree $d$ containing $X$. Denote by $Z$ the complete intersection of $Z_{1},\ldots,Z_{t}$. Then $Z$ is smooth.
\end{lem}

\begin{proof}
Let $\sigma$ be the linear system of hypersurfaces of
degree $d$ containing $X$. Then the base locus of $\sigma$ is just
$X$.

\medskip

If $t=1$, then $\text{Sing}(Z_{1})\subseteq X$ by
the Bertini's theorem. We can choose other $c-1$ hypersurfaces $Z_{1,
2},\ldots,Z_{1, c}$ of degree $d$ containing $X$ such that $X$ is
the complete intersection of $Z_{1},Z_{1, 2},\ldots,Z_{1, c}$. Then
Lemma \ref{lem: X cap Sing(Z) contained in Sing(X)} implies $\text{Sing}(Z_{1})\subseteq \text{Sing}(X)=\emptyset$, i.e. $Z_{1}$ is smooth.

\medskip

Assume $t\geq 2$ and the conclusion holds for $t-1$. Let $W$ be the complete
intersection of $Z_{1},\ldots,Z_{t-1}$. Then $W$ is smooth by induction. The restriction on $W$ of $\sigma$ gives a
linear system $\sigma'$ with base locus $X$ too. Since $Z_{t}$ is
general, $\text{Sing}(Z_{t}\cap_{sch} W)\subseteq X$ by the Bertini's theorem,
i.e. $\text{Sing}(Z)\subseteq X$. By Lemma \ref{lem: X cap Sing(Z) contained in Sing(X)}, $\text{Sing}(Z)\subseteq \text{Sing}(X)=\emptyset$, i.e. $Z$ is smooth.
\end{proof}

\bigskip

Let $X\subseteq\mathbb{P}^{N}$ be a projective variety, and take a fixed point $x\in X$. Denote by $\mathfrak{L}_{x, m}(X)=\{L\in\mathbb{G}(m, N)|x\in L\subseteq X\}$, where $\mathbb{G}(m, N)$ is the Grassmann variety parameterizing $m$-dimensional linear subspaces in $\mathbb{P}^{N}$. If $x\in X$ is a smooth point, then let $T_{x}X$ be the embedded affine tangent space of $X$ at $x$. Denote by $\mathbb{T}_{x}X$ the closure of $T_{x}X$ in
$\mathbb{P}^{N}$.

\begin{lem} \label{lem: T_X cap X covered by linear subspaces}
If $X\subseteq\mathbb{P}^{n+1}$ is a hyperquadric with $\text{dim}(X)=n\geq
2$, and $m$ is a positive integer such that $m\leq [\frac{n}{2}]$, then
for any point $x\in X$, we have $\mathbb{T}_{x}X\cap
X=\bigcup\limits_{L\in \mathfrak{L}_{x, m}(X)}L$.
\end{lem}

\begin{proof}
Induction on $m$. The case $m=1$ is trivial.

Assume $m>1$, then $n\geq 4$. If $x\in X_{sm}$, then
$\mathbb{T}_{x}X\cap X$ is a cone over an $(n-2)$-dimensional quadric $X'$ with the vertex $x$ (see for example, \cite[Prop. 2.2]{Ru}). By induction, for any
point $x'\in X'$, $\mathbb{T}_{x'}X'\cap X'=\bigcup\limits_{L'\in
\mathfrak{L}_{x', m-1}(X')}L'$. In particular, $\mathfrak{L}_{x',
m-1}(X')\neq\emptyset$. Take $L'\in\mathfrak{L}_{x', m-1}(X')$. Let
$L=\langle x, L'\rangle$ be the linear subspace spanned by $x$ and
$L'$, then $x'\in L\in\mathfrak{L}_{x, m}(X)$. Hence, $\mathbb{T}_{x}X\cap X=\bigcup\limits_{L\in \mathfrak{L}_{x,
m}(X)}L$.

\medskip

If $x\in \text{Sing}(X)$, then $X$ is a cone with $x$ being a vertex and
$\mathbb{T}_{x}X\cap X=X$. Take a general hyperplane section
$\widetilde{X}$ with $x\notin \widetilde{X}$. We can deduce that
$\mathbb{T}_{x}X\cap X=\bigcup\limits_{L\in \mathfrak{L}_{x,
m}(X)}L$ by taking the same argument as above where $X'$ is replaced
by $\widetilde{X}$.
\end{proof}

\begin{prop}[\cite{Sa}, Main Thm.] \label{prop: classification X swept out by large dimensional linear subspaces}
Let $X\subseteq\mathbb{P}^{N}$ be a smooth
projective variety of dimension $n\geq 2$.  Assume that $m$ is an integer such that
$m\geq \frac{n}{2}$. If for each point $x\in X$, there is an
$m$-dimensional linear subspace $L_{x}$ satisfying $x\in L_{x}\subseteq
X$, then $X$ is one of the following three cases:

$(i)$ a $\mathbb{P}^{n-d}$-bundle over a $d$-dimensional smooth
projective variety $S$ where a general linear subspace $L_{x}$ is in
the fiber of the canonical projection and $\frac{n}{2}\leq d\leq n$;

$(ii)$ an even-dimensional smooth hyperquadric;

$(iii)$ the Grassmann variety $\mathbb{G}(1, m+1)$
parameterizing lines in $\mathbb{P}^{m+1}$.
\end{prop}

\begin{proof}[Proof of Proposition \ref{prop: complete intersection of two hyperquadrics}]
Let $S$ be the homogeneous
coordinate ring of $\mathbb{P}^{n+2}$ and $I(X_{1})=(f_{1}),
I(X_{2})=(f_{2})$. By Lemma \ref{lem: general Z contain smooth X is smooth}, we can assume that $X_{1}$ and
$X_{2}$ are both smooth.

\medskip

$(i)$ It is just a special case of Theorem \ref{thm: introduction deg(X)|deg(Y)}.

\bigskip

$(ii)$ Assume $n$ to be even. $f_{1}, f_{2}$ are quadratic forms and
they uniquely correspond to symmetric coefficient matrices $A_{1}$
and $A_{2}$, respectively. The equation
$\text{det}(\lambda_{1}A_{1}+\lambda_{2}A_{2})=0$ has a solution for some
$[\lambda_{1}, \lambda_{2}]\in\mathbb{P}^{1}$. Equivalently, there
is a singular hyperquadric $X_{3}$ containing $X$. Thus, $X_{i}\neq X_{3}$, and $X=X_{i}\cap_{sch}X_{3}$ for $i=1, 2$. By Lemma \ref{lem: X cap Sing(Z) contained in Sing(X)}, $X\subseteq(X_{3})_{sm}$ and $X\cap\text{Sing}(X_{3})=\emptyset$. In particular, $\text{dim}(\text{Sing}(X_{3}))\leq 1$. However,  $\mathcal{O}_{X_{3}}(X)=l^{*}\mathcal{O}_{\mathbb{P}^{n+2}}(X_{1})$ is an ample Cartier divisor on $X_{3}$, where $l: X_{3}\rightarrow\mathbb{P}^{n+2}$ is the natural inclusion. Thus, $\text{dim}(\text{Sing}(X_{3}))=0$. So $X_{3}$ has only one singular point $p$, and $p\notin X$.

\medskip

Take any $x\in X$, then $x\in(X_{3})_{sm}$. Let $X_{3}'=H\cap X_{3}$, where $H$ is a general hyperplane of $\mathbb{P}^{n+2}$ passing through $x$. Then $X_{3}'$ is smooth. By Lemma \ref{lem: T_X cap X covered by linear subspaces}, there exists an $\frac{n}{2}$-dimensional subspace $L'$ of $H$ such that $x\in L'\subseteq X_{3}'$. Let $L=\langle p, L'\rangle$ be the linear subspace spanned by $p$ and $L'$ in $\mathbb{P}^{n+2}$. Then $L\subseteq X_{3}$ and $x\in L\cap_{sch}X_{1}$. Moreover, $L\subseteq X_{1}$, or $\text{dim}(L\cap_{sch}
X_{1})=\frac{n}{2}$.

\medskip

We claim that if $x\in X$ is general, then $L\nsubseteq X_{1}$, $L\cap_{sch} X_{1}$ is irreducible, and it has multiplicity one on the unique irreducible component. Otherwise, for a
general point $x\in X$, there is an $\frac{n}{2}$-dimensional linear
subspace $L_{x}$ satisfying that $x\in L_{x}\subseteq X$. On the
other hand,  $Chow_{\frac{n}{2}, 1}(X)$ is a projective scheme,
where $Chow_{\frac{n}{2}, 1}(X)$ is the Chow variety parameterizing
effective algebraic cycles of dimension $\frac{n}{2}$ and degree 1 on $X$. Hence, for an
arbitrary point $x\in X$, there is an $\frac{n}{2}$-dimensional
linear subspace $L_{x}$ satisfying that $x\in L_{x}\subseteq X$. By
Proposition \ref{prop: classification X swept out by large dimensional linear subspaces}, $X$ cannot be a smooth complete intersection of two
hyperquadrics.

\bigskip

$(iii)$ Fix an arbitrary point $x_{0}\in X$. Let
$\mathcal{Z}=\{V(f)|f=\lambda_{1}f_{1}+\lambda_{2}f_{2},\
[\lambda_{1}, \lambda_{2}]\in\mathbb{P}^{1}\}$. By Lemma \ref{lem: X cap Sing(Z) contained in Sing(X)}, for
any $Z\in\mathcal{Z}$, $x_{0}\in Z_{sm}$.

Let
$B=\overline{\bigcup\limits_{Z\in\mathcal{Z}}(\mathbb{T}_{x_{0}}Z\cap
X)}$. Note that for any $Z\in\mathcal{Z}$,
$\mathbb{T}_{x_{0}}Z\cap X$ is a divisor of $X$ and as a closed
subscheme of $\mathbb{P}^{n+2}$, $\text{deg}(\mathbb{T}_{x_{0}}Z\cap X)\leq
4$. If $B\subsetneq X$, then it is a possibly reducible variety of
pure dimension $n-1$, and there is a divisor $D$ with $Supp(D)\subseteq B$ such
that for general $Z\in\mathcal{Z}$, $\mathbb{T}_{x_{0}}Z\cap X=D$. However, if $Z_{1}, Z_{2}\in\mathcal{Z}$ are different
hyperquadrics, then $X=Z_{1}\cap_{sch} Z_{2}$ and
$\mathbb{T}_{x_{0}}Z_{1}\cap\mathbb{T}_{x_{0}}Z_{2}=\mathbb{T}_{x_{0}}X$. Thus, $\text{Supp}(D)\subseteq\mathbb{T}_{x_{0}}X$. Since $n\geq 3$, $\text{dim}(D)>\frac{\text{dim}(X)}{2}$. By Theorem \ref{thm: introduction deg(X)|deg(Y)}$(i)$, $\text{dim}\langle\text{Supp}(D)\rangle\geq n+1$. It is a contradiction. Hence, $X=B=\overline{\bigcup\limits_{Z\in\mathcal{Z}}(\mathbb{T}_{x_{0}}Z\cap
X)}\subseteq\overline{\bigcup\limits_{Z\in\mathcal{Z}}(\mathbb{T}_{x_{0}}Z\cap Z)}$.

\medskip

Consequently, for a general point $x\in X$, there exists some
$Z_{x}\in\mathcal{Z}$ such that $x\in\mathbb{T}_{x_{0}}Z_{x}\cap
Z_{x}$. By Lemma \ref{lem: T_X cap X covered by linear subspaces}, there exists an $[\frac{n+1}{2}]$-dimensional
linear subspace $L_{x}$ in $\mathbb{P}^{n+2}$ such that $x_{0}, x\in
L_{x}\subseteq Z_{x}$. Without loss of generality, we can assume that $Z_{x}\neq X_{1}$. By taking a similar argument with the last paragraph of the proof of $(ii)$, we can deduce from the generality of $x$ that $L_{x}\nsubseteq X$, $L_{x}\cap_{sch}X_{1}$ is irreducible, and it has multiplicity one on the unique irreducible component. Thus, $Q_{x}=L_{x}\cap_{sch}X_{1}\subseteq X$ is an $[\frac{n-1}{2}]$-dimensional quadric passing through $x_{0}$ and $x$.
\end{proof}

\section{\large Application to quadratic varieties} \label{section: applications to quadratic varieties}

Let $X\subsetneq\mathbb{P}^{N}$ be an $n$-dimensional nondegenerate smooth projective variety swept out by $m$-dimensional quadrics passing through a fixed point $x_{0}$. In \cite{KS}, Kachi and Sato showed  that if $m\geq\frac{3n}{5}+1, n=5,6,10$, or $m\geq\frac{3n}{5}, n\geq 4,
n\neq 5,6,10$, and each $Q_{x}$ is smooth, then $X$ is a hyperquadric. Then in \cite[Thm. 2]{Fu} it has been proved that if $m\geq[\frac{n}{2}]+2$, then $X\subseteq\mathbb{P}^{n+1}$ is a hyperquadric. Moreover in \cite[Prop. 3]{Fu} it has been shown that if $m\geq[\frac{n}{2}]+1$ with $n\geq 3$, and if $N\geq\frac{3n}{2}$, then $X\subseteq\mathbb{P}^{N}$ is projectively equivalent to one of the cases $(a)-(e)$ in Theorem \ref{thm: introduction application to var. swept out by large dim. quadrics}. Thus our Theorem \ref{thm: introduction application to var. swept out by large dim. quadrics} can be considered as a refinement of \cite[Prop. 3]{Fu}. For further generalizations of \cite{KS} and of \cite{Fu}, one can consult \cite{BI}.

\begin{proof}[Proof of Theorem \ref{thm: introduction application to var. swept out by large dim. quadrics}]
By \cite[Prop. 3]{Fu}, we can assume $N < \frac{3n}{2}$. If the Hartshorne Conjecture holds or if $X$ is
quadratic, then $X$ is
 a complete intersection in $\mathbb{P}^{N}$ by \cite{IR}. By Theorem \ref{thm: introduction deg(X)|deg(Y)}, we have $\text{deg}(X)=1, \text{ or } 2$. Since $X\subsetneq\mathbb{P}^{N}$ is nondegenerate, $\text{deg}(X)=2$ and $\text{dim}(L)=\text{dim}(X)+1$, where $L=\langle X\rangle$. Hence, $L=\mathbb{P}^{N}$ and $X$ is a hyperquadric.

\end{proof}

\small

INSTITUTE OF MATHEMATICS, AMSS, CHINESE ACADEMY OF SCIENCES,

\smallskip

55 ZHONGGUANCUN EAST ROAD, BEIJING, 100190, P. R. CHINA

\smallskip

E-mail address: qifengli@amss.ac.cn

\end{document}